\documentstyle[12pt,amssym]{article}

\newtheorem{theo}{Theorem}

\newtheorem{prop}[theo]{Proposition}
\newtheorem{coro}[theo]{Corollary}

\newcommand{\myC}{{\cal C}}

\newcommand{\myA}{{\cal A}}
\newcommand{\myB}{{\cal B}}

\newcommand{\bQ}{{\Bbb Q}}

\newcommand{\bR}{{\Bbb R}}
\newcommand{\bN}{{\Bbb N}}

\newcommand{\id}{{\rm id}}

\newcommand{\supp}{{\rm supp}}
\newcommand{\levy}{L\'{e}vy\,}

\begin{document}

\begin{center}
{\large\bf Convolution of probability measures on Lie groups and homogenous spaces}

Ming Liao\footnote{Department of Mathematics and Statistics, Auburn University, Auburn, AL 36849, USA.

Email: liaomin@auburn.edu}

\end{center}

\begin{quote}
{\bf Summary} \ We study (weakly) continuous convolution semigroups of probability measures on a Lie group $G$ or
a homogeneous space $G/K$, where $K$ is a compact subgroup. We show that such a convolution semigroup is the convolution
product of its initial measure and a continuous convolution semigroup with initial measure at the identity of $G$ or the origin of $G/K$.
We will also obtain an extension of Dani-McCrudden's result on embedding an infinitely divisible probability measure in
a continuous convolution semigroup on a Lie group to a homogeneous space.



\end{quote}

\section{Introduction} \label{sec1}

Let $G$ be a Lie group. For a measure $\mu$ and a function $f$, the integral $\int fd\mu$, whenever it is defined, may be written as $\mu(f)$. For any two measures $\mu$ and $\nu$ on $G$ (defined on Borel sets), their convolution $\mu*\nu$
is the measure on $G$ defined as usual by
\begin{equation}
\mu*\nu(f) = \int_{G\times G} f(xy)\mu(dx)\nu(dy) \label{mu*nuG}
\end{equation}
for any $f\in\myB(G)_+$, where $\myB(G)_+$ is the space of nonnegative Borel functions on $G$.

  The convolution may be defined on the homogeneous space $X=G/K$, where $K$ is a compact subgroup.
In \cite{gangolli}, the convolution on $X$ is defined by identifying the measures on $X$ as $K$-right invariant measures
on $G$ (see Proposition~\ref{prmunu}), and use the convolution on $G$. We give an equivalent but more direct definition
here. Let $\pi$: $G\to G/K$ be the natural projection. A Borel measurable map $S$: $X\to G$ is called a section map on $X$
if $\pi\circ S=\id_X$, the
identity map on $X$. Such an $S$ is not unique and may not be continuous, but it can always be chosen to be smooth near any
point in $X$. The convolution $\mu*\nu$ of any two measures $\mu$ and $\nu$ on $X$ may be defined by
\begin{equation}
\mu*\nu(f) = \int_{X\times X}\int_Kf(S(x)ky)\rho_K(dk)\mu(dx)\nu(dy), \label{mu*nuX}
\end{equation}
for any $f\in\myB(X)_+$, where $\rho_K$ is the normalized Haar measure on $K$. This definition is independent
of the choice for the section map $S$, and reduces to the convolution on $G$ when $K=\{e\}$, where $e$ is the identity
element of $G$. Because for any $g\in G$ and $x\in X$, $S(gx)=gS(x)k$ for some $k\in K$, it is
easy to show that the convolution on $X$, as on $G$, is associative, that
is, $(\mu*\nu)*\gamma=\mu*(\nu*\gamma)$. Therefore, the $n$-fold convolution product $\mu^{*n}=\mu*\mu*\cdots*\mu$
is well defined.

  Convolution of functions on $X=G/K$ have appeared in literature under various contexts, but in the present setting, they
all take the form below (see for example \cite{vv}): For $f_1,f_2\in\myB(X)_+$,
\begin{equation}
f_1*f_2(gK) = \int_G f_1(hK)f_2(h^{-1}gK)\rho(dh), \quad g\in G, \label{f1*f2X}
\end{equation}
where $\rho$ denotes a left invariant Haar measure on $G$. It is easy to show that this definition
is compatible with our definition of convolution of measures in the sense that if $\mu_1$ and $\mu_2$ are measures on $X$
with densities $f_1$ and $f_2$ with respect to $\pi\rho$, then $\mu_1*\mu_2$ has density $f_1*f_2$, where $\pi\rho$
is the measure on $X$ defined by $\pi\rho(f)=\rho(f\circ\pi)$ for $f\in\myB(X)_+$.
%% To prove this, note that any left invariant Haar measure $\rho$ on $G$ is $K$-right invariant. Indeed,
% for $k\in K$, $r_k\rho$ is left invariant and so $r_k\rho=\lambda(k)\rho$ for some $\lambda(k)\in\bR$. For
% any continuous function $f$ with compact support on $G$, $\rho(f\circ r_k)=\lambda(k)\rho(f)$. This shows $k\mapsto
% \lambda(k)$ is continuous. Then $k\mapsto\lambda(k^{-1})$ is a continuous Lie group homomorphism from $K$ into the
% multiplicative group $(0,\,\infty)$. Its range as a compact subgroup is necessarily $\{1\}$, and hence $r_k\rho=\rho$.

A family of probability measures $\mu_t$ on $G$ or on $X$, $t\in\bR_+=[0,\,\infty)$, is called a convolution semigroup
on $G$ or on $X$ if $\mu_{t+s}=\mu_t*\mu_s$ for $s,t\in\bR_+$. It is called continuous if $\mu_t\to\mu_0$ weakly
as $t\to 0$. This in fact implies $\mu_t\to\mu_s$ weakly as $t\to s$ for any $s>0$ (see \cite[Theorem~1.5.7]{heyer}
for $\mu_t$ on $G$, then the result on $X$ follows from Proposition~\ref{prmutnut}).

  Convolution semigroups on Lie groups or on more general topological groups have been studied extensively, often
in connection with parabolic equations or Markov processes. We will not attempt to provide a review of literature
on this subject, except to mention that a systematic study of convolution semigroups on locally compact groups may
be found in Heyer \cite{heyer}, and some of more recent developments in Lie group setting in Applebaum \cite{apple2}.

  Most studies are centered on continuous convolution semigroups $\mu_t$ on $G$ with $\mu_0=\delta_e$, the
unit point mass at $e$. In this case, we may associate a Markov process $g_t$ in $G$ with distribution $\mu_t$.
By shifting the process by left translations, we obtain a family of processes, one for each starting point in $G$,
and all governed by the transition semigroup of operators $P_t$, defined by $P_tf(g)=\int f(gh)\mu_t(dh)$
for $f\in\myB(G)_+$ and $t\geq 0$. The $P_t$ is in fact a so-called Feller transition semigroup, and as a Feller
process, $g_t$ has many useful properties, such as the strong Markov property, and after a suitable modification on
a null set for each $t>0$, it will have rcll paths (right continuous paths with left limits). Moreover, for
any $s<t$, $g_s^{-1}g_t$ is independent of the process up to time $s$ and has the distribution $\mu_{t-s}$. Such a process
is called a \levy process in $G$, see \cite{Liao} for more details.

When $\mu_t$ is a convolution semigroup on $G$, but $\mu_0\neq\delta_e$, we may still associate a process, but it
is not clear if this is a Feller process or can be started at an arbitrary point in $G$.

A natural question is whether any continuous convolution semigroup $\mu_t$ on $G$ can be expressed as $\mu_t=\mu_0*\mu_t^e$
for a continuous convolution semigroup $\mu_t^e$
with $\mu_0^e=\delta_e$. This would allow us to associate a Feller process with distribution $\mu_t$,
whose transition semigroup $P_t$ is determined by $\mu_t^e$.
We have not been able to find a related discussion in the literature.

A purpose of this paper is to give an affirmative answer to this question in section \S2.
We will also show in \S2 that if $\mu_t$ is a continuous convolution semigroup on $X=G/K$, then there is
a continuous convolution semigroup $\mu_t^o$ with $\mu_0^o=\delta_o$, the unit point mass at the origin $o=eK$
of $G/K$, such that $\mu_t=\mu_0*\mu_t^o$.

A probability measure is also called a distribution. A distribution $\mu$ on $G$ or on $X$ is called infinitely divisible
if for any $n\in\bN$ (the set of natural numbers), there is a distribution $\nu$, called an $n$th root of $\mu$,
with $\mu=\nu^{*n}$. If $\mu_t$ is a convolution semigroup on $G$ or $X$, then $\mu=\mu_1$
is clearly infinitely divisible. It is a challenging problem to determine whether
any infinitely divisible distribution $\mu$ can be embedded in a continuous convolution semigroup, that is,
whether there is a continuous convolution semigroup $\mu_t$ with $\mu=\mu_1$.

  The embedding problem on locally compact groups has been studied extensively, see Heyer \cite[chapter~III]{heyer}
for the results up to 1976. More recently, the problem was solved by Dani-McCrudden \cite{dm0,dm} for
a large class of Lie groups. They call a Lie group $G$ to belong to class $\myC$ if $G$ is connected and there is
a linear action of $G$ on a finite dim real vector space $V$ with a discrete kernel, that is, the set of $g\in G$ that
acts as the identity map on $V$ is a discrete subset of $G$. Note that a connected compact or semi-simple Lie
group is of class $\myC$.

In \S3, we will show that an infinitely divisible distribution on $G/K$ can be embedded
in a continuous convolution semigroup on $G/K$ if $G$ is of class $\myC$, thus extending Dani-McCrudden's result
from Lie groups to homogeneous spaces. This allows us to reduce the study of an infinitely divisible distribution to that
of a continuous convolution semigroup which is more tractable. For example, Gangolli \cite{gangolli} obtained
an interesting \levy-Khinchin type formula for infinitely divisible distributions on symmetric spaces,
in terms of spherical functions. The same formula was obtained by Applebaum \cite{apple}
and Liao-Wang \cite{LiaoWang}
for convolution semigroups by a much simpler method based on Hunt's generator formula \cite{hunt}.
Because a symmetric space is a homogeneous space of a semi-simple Lie group, to which the embedding result holds,
we may now derive the more difficult result in \cite{gangolli} by the easier method in \cite{apple,LiaoWang}.
% and \cite{LiaoWang}.

  The embedding on $G/K$ is equivalent to embedding an infinitely divisible distribution on $G$, that
is $K$-right invariant and has a $K$-right invariant $n$th root for any $n$, in
a continuous convolution semigroup on $G$ that is also $K$-right invariant. Some feature of the embedding on $G$
given in \cite{dm} can be used to show such an embedding is possible.

We note that the embedding may hold on Lie groups that are not of class $\myC$, even some locally compact
groups, see for example \cite{dgs} for some of such Lie groups. It may be possible to
obtain the embedding on the homogeneous spaces of these groups by the present approach.

\section{Continuous convolution semigroups} \label{sec2}

For $g\in G$, let $l_g$, $r_g$ and $c_g$ be, respectively, the left translation $x\mapsto gx$, the
right translation $x\mapsto xg$, and the conjugation map $x\mapsto gxg^{-1}$ on $G$. A measure $\mu$ on $G$ is called left invariant if $l_g\mu=\mu$ for any $g\in G$. It is called $K$-left invariant if $l_k\mu=\mu$ for any $k\in K$. Similarly,
we may define right invariant, $K$-right invariant, conjugate invariant and $K$-conjugate invariant measures on $G$.
A measure $\mu$ on $G$ is called bi-invariant (resp. $K$-bi-invariant) if it is both left and right (resp.
both $K$-left and $K$-right) invariant.

  The Lie group $G$ acts on $X=G/K$ naturally by $(g,xK)\mapsto gxK$ for $g,x\in G$. A measure $\mu$ on $X$
is called $K$-invariant if $k\mu=\mu$ for $k\in K$. It is easy to see that if $\nu$ is $K$-invariant, then the convolution $\mu*\nu$, defined by (\ref{mu*nuX}), may be written more concisely as
\begin{equation}
\mu*\nu(f) = \int_{X\times X}f(S(x)y)\mu(dx)\nu(dy), \label{mu*nuX2}
\end{equation}
which is independent of the choice for the section map $S$. Moreover, if $\mu$ is also $K$-invariant, then so
is $\mu*\nu$.

\begin{prop} \label{prmunu}
  The map
\[\mu\ \mapsto\ \nu=\pi\mu\]
is a bijection from the set of $K$-right invariant measures $\mu$ on $G$ onto the set of measures $\nu$ on $X=G/K$.
It is also a bijection from the set of $K$-bi-invariant measures $\mu$ on $G$ onto the set of $K$-invariant measures $\nu$
on $X$. Moreover, if $\nu$ is a measure on $X$, then the unique $K$-right invariant measure $\mu$
on $G$ satisfying $\nu=\pi\mu$ is given by
\begin{equation}
\mu(f) = \int_X\int_K f(S(x)k)\rho_K(dk)\,\nu(dx) \quad \mbox{for $f\in\myB(G)_+$,} \label{15munu}
\end{equation}
where $S$ is any section map on $X$. Furthermore, the map $\mu\to\pi\mu$ preserves the convolution in the sense that
for measures $\mu_1$ and $\mu_2$ on $G$,
\begin{equation}
\pi(\mu_1*\mu_2) = (\pi\mu_1)*(\pi\mu_2), \label{pimu*pinu}
\end{equation}
provided one of the following three conditions holds: $\mu_1$ is $K$-right invariant, or $\mu_2$ is $K$-left invariant,
or $\mu_2$ is $K$-conjugate invariant.
\end{prop}

\noindent {\bf Proof} \ For $g\in G$, $S\circ\pi(g)=gk$ for some $k\in K$. If $\mu$ is a $K$-right invariant measure on $G$
with $\pi\mu=\nu$, then for $f\in\myB(G)_+$,
\[\mu(f) = \int_K\mu(f\circ r_k)\rho_K(dk) = \int_K\mu(f\circ r_k\circ S\circ\pi)\rho_K(dk) =
\int_K\nu(f\circ r_k\circ S)\rho_K(dk).\]
This shows that $\mu$ satisfies (\ref{15munu}). Conversely, using the $K$-right invariance of $\rho_K$, it is
easy to show that $\mu$ given by (\ref{15munu}) is $K$-right invariant with $\pi\mu=\nu$. Because
for $g,h\in G$, $gS(h)=S(gh)k$ for some $k\in K$, and $g\pi(h)=\pi(gh)$, it is easy to show that $\mu$
is $K$-bi-invariant if and only if $\nu$ is $K$-invariant. For measures $\mu_1$ and $\mu_2$ on $G$, satisfying one
of the three conditions stated in the proposition, let $\nu_1=\pi\mu_1$ and $\nu_2=\pi\mu_2$. Then for $f\in\myB(X)_+$,
\begin{eqnarray*}
&& \pi(\mu_1*\mu_2)(f) = \int f(\pi(g_1g_2))\mu_1(dg_1)\mu_2(dg_2) = \int f(\pi(g_1kg_2))\mu_1(dg_1)\rho_K(dk)\mu_2(dg_2) \\
&=& \int f(g_1k\pi(g_2))\mu_1(dg_1)\rho_K(dk)\mu_2(dg_2) = \int f(gky)\mu_1(dg)\rho_K(dk)\nu_2(dy) \\
&=& \int f(S(\pi(g))k'ky)\mu_1(dg)\rho_K(dk)\nu_2(dy) \quad \mbox{(for some $k'\in K$)} \\
&=& \int f(S(x)ky)\nu_1(dx)\rho_K(dk)\,\nu_2(dy) = \nu_1*\nu_2(f). \ \ \Box
\end{eqnarray*}

For a convolution semigroup $\mu_t$ on $G$, because $\mu_0*\mu_0=\mu_0$, by \cite[Theorem~1.2.10]{heyer}, $\mu_0=\rho_H$
for some compact subgroup $H$ of $G$. Then $\mu_t$ is $H$-bi-invariant because $\mu_t=\mu_0*\mu_t=\mu_t*\mu_0$.
It now follows from Proposition~\ref{prmunu} that
if $\nu_t$ is a convolution semigroup on $X=G/K$, then $\nu_0=\pi\rho_H$ for some compact subgroup $H$ of $G$ containing $K$.

\begin{prop} \label{prmutnut}
(a) \ If a convolution semigroup $\mu_t$ on $G$ is $K$-right invariant, that is, each $\mu_t$ is $K$-right invariant, then
it is $K$-bi-invariant.

\noindent (b) \ Any convolution semigroup $\nu_t$ on $X=G/K$ is $K$-invariant.

\noindent (c) \ The map
\[\mu_t\ \mapsto\ \nu_t=\pi\mu_t\]
is a bijection from the set of $K$-bi-invariant convolution semigroups $\mu_t$ on $G$ onto the set
of convolution semigroups $\nu_t$ on $X=G/K$. Moreover, $\mu_t$ is continuous if and only if so is $\nu_t$.
\end{prop}

\noindent {\bf Proof:} \ By the preceding discussion, $\mu_0=\rho_H$ and $\mu_t$ is $H$-bi-invariant
for some compact subgroup $H$ of $G$. The $K$-right invariance of $\mu_0$ implies $K\subset H$. This proves (a). Let $\nu_t$ be a convolution semigroup on $X$, and let $\mu_t$ be the unique $K$-right invariant probability measure on $G$
with $\pi\mu_t=\nu_t$. By Proposition~\ref{prmunu} and (a), $\mu_t$ is a $K$-bi-invariant convolution semigroup on $G$,
and hence $\nu_t=\pi\mu_t$ is $K$-invariant. This proves (b). Now (c) follows from (b) and Proposition~\ref{prmunu},
but note that to derive the continuity of $\mu_t$ from that of $\nu_t$, (\ref{15munu}) is used,
where $\int_Kf(S(x)k)\rho_K(dk)$ is continuous in $x$ if $f$ is bounded continuous on $G$. This is
because this integral is independent of $S$ and $S$ may be chosen to be continuous near any $x$. \ $\Box$
\vspace{2ex}

% Let $\mu$ be a $K$-conjugate invariant measure on $G$, that is, $c_k\mu=\mu$ for $k\in K$. Then $\nu=\pi\mu$
% is a $K$-invariant measure on $X$, because for any $f\in\myB(X)_+$ and $k\in K$,
% \[\nu(f\circ k) = \mu(f\circ k\circ\pi) = \mu(f\circ\pi\circ l_k) = \mu(f\circ\pi\circ c_k) =
% \mu(f\circ\pi) = \nu(f).\]
%
\begin{prop} \label{prmu2nu}
Let $\mu_t$ be a continuous $K$-conjugate invariant convolution semigroup on $G$ (that is, each $\mu_t$ is $K$-conjugate
invariant). Then $\nu_t=\pi\mu_t$ is a continuous convolution semigroup on $X$. Note that
if $\mu_0=\delta_e$, then $\nu_0=\delta_o$.
\end{prop}

\noindent {\bf Proof:} \ We just have to show $\nu_{s+t}=\nu_s*\nu_t$, which follows from Proposition~\ref{prmunu}. \ $\Box$
% For $f\in\myB(X)_+$,
% \begin{eqnarray*}
% && \nu_{s+t}(f) = \mu_{s+t}(f\circ\pi) = \mu_s*\mu_t(f\circ\pi) = \int(f\circ\pi)(gh)\mu_s(dg)\mu_s(dh) \\
% &=& \int f(g\pi(h))\mu_s(dg)\mu_t(dh) = \int f(gy)\mu_s(dg)\nu_t(dy) = \int f(S(\pi(g))ky)\mu_s(dg)\nu_t(dy) \\
% && \mbox{($g=S(\pi(g))k$ for any section map $S$, and some $k\in K$ depending on $g$ and $S$)} \\
% &=& \int f(S(\pi(g))y)\mu_s(dg)\nu_t(dy) \quad \mbox{(because $\nu_t$ is $K$-invariant)} \\
% &=& \int f(S(x)y)\nu_s(dx)\nu_t(dy) = \nu_s*\nu_t(f). \quad \Box
% \end{eqnarray*}
\vspace{2ex}

Recall that for any continuous convolution semigroup $\mu_t$ on $G$ with $\mu_0=\delta_e$, there is
a \levy process $g_t$ in $G$ with distribution $\mu_t$ and rcll paths, whose transition semigroup $P_t$ is given
by $P_tf(g)=\int f(gh)\mu_t(dh)$ for $f\in\myB(G)_+$ and $g\in G$. Note that $P_t$ is left invariant
in the sense that $P_t(f\circ l_g)=(P_tf)\circ l_g$ for $g\in G$. When $\mu_t$ is $K$-conjugate invariant, then $P_t$
is also $K$-right invariant in the sense that $P_t(f\circ r_k)=(P_tf)\circ r_k$ for $k\in K$.

Let $\nu_t$ be a continuous convolution semigroup on $X=G/K$ with $\nu_0=\delta_o$. We may associate
a Markov process $x_t$ in $X$ with distribution $\nu_t$ and rcll paths such that its transition semigroup $Q_t$
is $G$-invariant in the sense that $Q_t(f\circ g)=(Q_tf)\circ g$ for $f\in\myB(X)_+$ and $g\in G$,
where $Q_tf(x)=\int_Xf(S(x)y)\nu_t(dy)$, which is independent of the choice of the section map $S$. Such
a process $x_t$ is called a $G$-invariant Markov process in $X$, see
\cite[\S2.2]{Liao}.

By \cite[Theorem~2.2]{Liao}, given a $G$-invariant Markov process $x_t$ in $X$ with $x_0=o$, there
is a \levy process $g_t$ in $G$ with $g_0=e$, that is also $K$-right invariant, such that $x_t$ is equal to $\pi(g_t)$
in distribution. Because the distribution $\mu_t$ of $g_t$ is $K$-conjugate invariant, the above may be stated as
a result for convolution semigroups, providing a converse to Proposition~\ref{prmu2nu}, and is recorded as follows.

\begin{theo} \label{thnu2mu}
Let $\nu_t$ be a continuous convolution semigroup on $X=G/K$ with $\nu_0=\delta_o$. Then there is
a continuous $K$-conjugate invariant convolution semigroup $\mu_t$ on $G$ with $\mu_0=\delta_e$ such that $\nu_t=\pi\mu_t$.
\end{theo}

Let $\mu_t$ be a continuous $K$-conjugate invariant convolution semigroup on $G$ with $\mu_0=\delta_e$. Then it is
easy to show that $\mu_t*\rho_K=\rho_K*\mu_t$, and $\mu_t'=\rho_K*\mu_t$ is a continuous $K$-bi-invariant convolution semigroup on $G$. By the result below, any continuous convolution semigroup on $G$ can be obtained this way, thus giving
an affirmative answer to the question raised in \S1.

\begin{theo} \label{thrhomute}
Let $\mu_t$ be a continuous convolution semigroup on $G$. Then $\mu_0=\rho_K$ for some compact subgroup $K$ of $G$,
and $\mu_t=\rho_K*\mu_t^e=\mu_t^e*\rho_K$ for some continuous $K$-conjugate invariant convolution semigroup $\mu_t^e$
with $\mu_0^e=\delta_e$. Consequently, $\mu_t$ is $K$-bi-invariant.
\end{theo}

\noindent {\bf Proof:} \ We have seen that $\mu_0=\rho_K$ for some compact subgroup $K$
of $G$. By the $K$-bi-invariance of $\rho_K$ and $\mu_t=\mu_0*\mu_t=\mu_t*\mu_0$, $\mu_t$ is $K$-bi-invariant.
By Proposition~\ref{prmutnut}, $\nu_t=\pi\mu_t$ is a continuous convolution semigroup on $X=G/K$
with $\nu_0=\pi\mu_0=\delta_o$. By Theorem~\ref{thnu2mu}, there is a continuous $K$-conjugate
invariant convolution semigroup $\mu_t^e$ on $G$ with $\mu_0^e=\delta_e$ such that $\pi\mu_t^e=\nu_t$. Then $\mu_t'= \rho_K*\mu_t^e=\mu_t^e*\rho_K$ is a continuous $K$-bi-invariant convolution semigroup on $G$ with $\pi(\mu_t')=\nu_t$.
By the uniqueness in Proposition~\ref{prmutnut}\,(c), $\mu_t'=\mu_t$. \ $\Box$
\vspace{2ex}

\noindent {\bf Remark:} \ We note that $\mu_t^e$ in Theorem~\ref{thrhomute} is not unique. This is because $\mu_t$
in Theorem~\ref{thnu2mu} is not unique. By \cite[Proposition~2.7]{Liao}, there is a family of $K$-conjugate
invariant \levy processes $g_t$ in $G$, all starting at $e$ but with different distributions, such that the
processes $x_t=\pi(g_t)$ in $X$ have the same distribution. The non-uniqueness of $\mu_t^e$ reflects the following fact:
For a continuous convolution semigroup $\mu_t$ on $G$ with $\mu_0=\rho_K$, where $K\neq\{e\}$, there are different \levy processes $g_t$ in $G$, different in the sense that they have different transition semigroup $P_t$, such that if they
all start with the same initial distribution $\rho_K$, then their $1$-dimensional distributions are the same, and are given
by $\mu_t$.
\vspace{2ex}

For two compact subgroups $K\subset H$ of $G$, the homogeneous space $H/K$ is a closed topological subspace
of $X=G/K$ (\cite[II.Proposition~4.4]{helgason1}), and $\pi\rho_H$ is the unique $H$-invariant
probability measure on $G/K$ supported by $H/K$, where $\pi$ is the natural projection: $G\to G/K$ as before. The
result below is an extension of Theorem~\ref{thrhomute} from $G$ to $G/K$.

\begin{theo} \label{thrhomuto}
Let $\nu_t$ be a continuous convolution semigroup on $X=G/K$. Then $\nu_0=\pi\rho_H$ for a compact subgroup $H$
of $G$ containing $K$, and $\nu_t=\nu_0*\nu_t^o=\nu_t^o*\nu_0$ for a continuous convolution semigroup $\nu_t^o$ on $X$
with $\nu_0^o=\delta_o$. Consequently, $\nu_t$ is $H$-invariant.
\end{theo}

\noindent {\bf Proof} \ By Proposition~\ref{prmutnut}, there is a
unique $K$-bi-invariant continuous convolution semigroup $\mu_t$ on $G$ such that $\nu_t=\pi\mu_t$.
By Theorem~\ref{thrhomute}, $\mu_0=\rho_H$ for some compact subgroup $H$ of $G$, and $\mu_t= \rho_H*\mu_t^e=\mu_t^e*\rho_H$
is bi-$H$-invariant, where $\mu_t^e$ is a continuous $H$-conjugate
invariant convolution semigroup on $G$ with $\mu_0^e=\delta_e$. The $K$-bi-invariance of $\rho_H=\mu_0$
implies $K\subset H$. Let $\nu_t^o=\pi\mu_t^e$. Then $\nu_t^o$ is a continuous convolution semigroup
on $X$ with $\nu_0^o=\delta_o$, and $\nu_t=\pi\mu_t=\nu_0*\nu_t^o=\nu_t^o*\nu_0$. \ $\Box$

\section{Embedding problem} \label{sec3}

As mentioned in \S1, by \cite{dm0,dm}, any infinitely divisible distribution on a Lie group $G$ of class $\myC$ can
be embedded in a continuous convolution semigroup on $G$. The following result provides its extension to homogeneous spaces.

\begin{theo} \label{thembeddingX}
Let $G$ be a Lie group of class $\myC$ and let $K$ be a compact subgroup. Then any infinitely divisible distribution $\alpha$
on $X=G/K$ can be embedded in a continuous convolution semigroup $\alpha_t$ on $X$, that is, $\alpha=\alpha_1$.
\end{theo}

\noindent {\bf Proof:} \ First let $\mu$ be an infinitely divisible distribution on $G$. Choose an $n$th root
of $\mu$, denoted as $r(n)$, for each $n\in\bN$. Then $n\mapsto r(n)$ is a map from $\bN$ to the set of roots
of $\mu$, called a root map of $\mu$. Let $Z^0(\mu,[G,G])$ be the subgroup of $G$ defined in \cite{dm},
whose precise definition is not important here,
we just need to know that it is a set of points $z\in G$ that centralize points in $\supp(\mu)$, the support
of $\mu$ (defined to be the smallest closed set $F$ with $\mu(F^c)=0$), that is, $zxz^{-1}=x$ for any $x\in\supp(\mu)$.
Let $\Gamma$ be the closure of convolution products of the measures of the following form, under the weak convergence topology,
\[zr(m)z^{-1} = (l_z\circ r_{z^{-1}})r(m); \quad z\in Z^0(\mu,[G,G])\ {\rm and}\ m\in\bN.\]
It is stated in \cite[Theorem~1.3]{dm} that, for some $n\in\bN$, there is an $n$th root $\nu$ of $\mu$ such that $\nu$ is embeddable in a convolution semigroup $\nu_t$ on $G$ with $\nu_t
\in\Gamma$. In fact, the proof of this theorem shows that $\nu$ can be embedded in
a convolution semigroup on rational times, and then the reader is
referred to \cite{heyer} for additional details to complete the proof.

Because this feature of the embedding in \cite{dm} is important for our purpose,
for the convenience of a reader who is not an expert in the embedding problem,
we will provide some additional details, according to the directions in the last
paragraph in \cite{dm}, to show that there is a continuous convolution semigroup $\nu_t$
with $\nu_1=\nu$ such that each $\nu_t$ is contained in $\tilde{\Gamma}$, the closure
of the set of measures of the form $x\gamma$ with $x\in G$ and $\gamma\in\Gamma$. It then follows that $\mu$, which has $\nu$ as an $n$th root, can be embedded
in the continuous convolution semigroup $\mu_t=\nu_t^{*n}$. We will do this later,
but for now, will proceed to the proof of Theorem~\ref{thembeddingX}.

Let $\alpha$ be an infinitely divisible distribution on $X=G/K$. For any $m\in\bN$, let $\alpha'$ be an $m$th
root of $\alpha$. Let $\mu$ and $\mu'$ be respectively the unique $K$-right invariant distributions on $G$
with $\alpha=\pi\mu$ and $\alpha'=\pi\mu'$ (see Proposition~\ref{prmunu}), and let a root map $r$ of $\mu$ be defined
by $r(m)=\mu'$. By the preceding discussion, there is a continuous convolution semigroup $\mu_t$ on $G$
with $\mu=\mu_1$ and $\mu_t=\nu_t^{*n}$. Because $\mu$ is $K$-right invariant, for any $x\in\supp(\mu)$
and $k\in K$, $xk\in\supp(\mu)$. Because $z\in Z^0(\mu,[G,G])$ centralizes points in $\supp(\mu)$, $xzkz^{-1}=
z(xk)z^{-1}=xk$, and so $zkz^{-1}=k$. Thus, $z$ also centralizes $k$. Because $r(m)$ is $K$-right invariant, it is now
easy to derive the $K$-right invariance
of $zr(m)z^{-1}$ from $k\in K$ being centralized by $z$. It follows that any measure in $\Gamma$
is $K$-right invariant, so is any measure in $\tilde{\Gamma}$. Because $\nu_t\in\tilde{\Gamma}$, $\mu_t=\nu_t^{*n}$
is $K$-right invariant. By Proposition~\ref{prmunu}, $\alpha_t=\pi\mu_t$ is a convolution semigroup on $X$,
and is continuous because so is $\mu_t$. We have $\alpha=\pi\mu=\pi\mu_1=\alpha_1$,
and this proves Theorem~\ref{thembeddingX}.

We now provide some details for embedding $\nu$ in a continuous convolution semigroup $\nu_t$ on $G$
with $\nu_t\in\tilde{\Gamma}$ as mentioned above. Let $\bQ$ be the set of rational numbers.
A family of distributions $\gamma_t$ on $G$, $t\in\bQ_+^*=
\bQ\cap(0,\,\infty)$, such that $\gamma_t*\gamma_s=\gamma_{t+s}$ for $t,s\in\bQ_+^*$, is called a convolution semigroup
on $G$ parameterized by $\bQ_+^*$. A convolution semigroup
parameterized by $\bQ_+=\bQ\cap[0,\,\infty)$ is defined similarly. In this latter case, $\gamma_0*\gamma_0=\gamma_0$, and
hence, $\gamma_0=\rho_H$ for some compact subgroup $H$ of $G$. Such
a convolution semigroup $\gamma_t$ is called continuous if it is continuous in $t$ under the weak convergence.
For $\gamma_t$ parameterized by $\bQ_+$, the continuity is equivalent to its continuity at $t=0$ (\cite[Theorem~1.5.7]{heyer}).
A distribution $\gamma$ on $G$ is called $\bQ_+$- (resp. $\bQ_+^*$-) embeddable if there is
a convolution semigroup $\gamma_t$ parameterized by $\bQ_+$ (resp. by $\bQ_+^*$) with $\gamma=\gamma_1$. The $\gamma$ is called
root compact if the set of distributions $\xi^{*k}$, where $\xi$ is an $n$th root of $\gamma$ for some $n\in\bN$
and $1\leq k\leq n$, is relatively compact under the weak convergence. This definition
agrees with \cite{heyer}, but it is called strongly root compact in \cite{dm}.

It is shown at the end of the proof of \cite[Theorem~1.3]{dm} that the $\nu$ above is root compact, and there is
a convolution semigroup $\gamma_t$ on $G$ parameterized by $\bQ_+^*$ with $\gamma_t\in\Gamma$ and $\gamma_1=\nu$.

Let $\myA$ be the set of weak limits of $\gamma_t$ as $t\to 0$, $t\in\bQ_+^*$. We now show
\begin{equation}
\forall\beta\in\myA, \quad \gamma_t*\beta = \beta*\gamma_t \quad \mbox{for $t\in\bQ_+^*$.} \label{gammatbeta}
\end{equation}
Let $t_n\in\bQ_+^*$ with $t_n\to 0$ such that $\gamma_{t_n}\to\beta$ weakly. Because $\gamma_1$ is root compact,
for any fixed and finite $T>0$, the set $\{\gamma_t$; $t\in\bQ_+^*\cap(0,\,T]\}$ is relatively weakly compact.
For $t\in\bQ_+^*$, $\gamma_t=\gamma_{t_n}*\gamma_{t-t_n}=\gamma_{t-t_n}*\gamma_{t_n}$.
Along a subsequence, $\gamma_{t-t_n}$ converges to some distribution $\lambda$ on $G$, so $\gamma_t=
\beta*\lambda=\lambda*\beta$. Then $\beta*\gamma_t=\beta*\beta*\lambda=\beta*\lambda*\beta=\gamma_t*\beta$. This
proves (\ref{gammatbeta}). The proof also shows that any $\beta\in\myA$ is a left and right convolution factor of $\gamma_t$.

For any closed subgroup $H$ of $G$, let $N(H)=\{x\in G$; $xHx^{-1}=H\}$. This is a closed subgroup
of $G$, called the normalizer of $H$, which contains $H$ as a normal subgroup.
By \cite[Theorems 3.5.1 and 3.5.4]{heyer} and their proofs, $\myA$ is a compact and connected abelian group under the convolution product
and the weak convergence topology, with identity $\rho_H$ for some compact subgroup $H$ of $G$,
and $\myA=\{x\rho_H$; $x\in H'\}$,
where $H'$ is a compact subgroup with $H\subset H'\subset N(H)$. Moreover, $H'/H$ and $\myA$ are isomorphic as topological groups under the map: $xH\mapsto x\rho_H$ from $H'/H$ to $\myA$, and there is a continuous convolution semigroup $\nu_t'$
parameterized by $\bQ_+$ with $\nu_0'=\rho_H$ such that for $t\in\bQ_+^*$, $\nu_t'=x_t\gamma_t=\delta_{x_t}*\gamma_t$
for some $x_t\in H_0'$, where $H_0'$ is the identity component of $H'$.

  Note that $x\rho_H=\rho_H x$ for $x\in H'$. To show this, just observe that
because $x$ normalizes $H$, $x\rho_H x^{-1}$ is a probability measure supported by $H$ and is left $H$-invariant, so must be equal to $\rho_H$. For $x\in H'$,
because $x\rho_H\in\myA$, by (\ref{gammatbeta}), $(x\rho_H)*\gamma_t=\gamma_t*(x\rho_H)$. As shown in the proof of
(\ref{gammatbeta}), any $\beta\in\myA$ is a left and right convolution factor of $\gamma_t$, this implies, with $\beta=
\rho_H$, that $\gamma_t=\rho_H*\gamma_t=\gamma_t*\rho_H$. Then for $x\in H'$,
\[x\gamma_t = x\rho_H*\gamma_t = \gamma_t*(x\rho_H) = \gamma_t*(\rho_H x) = \gamma_tx.\]
Because $\myA$ is abelian, we obtain $(x\gamma_t)*(y\gamma_s) = xy\gamma_{s+t} = yx\gamma_{s+t}$ for $x,y\in H'$
and $s,t\in\bQ_+^*$.

It is well know that the exponential map is surjective on a compact connected Lie group, so there is a vector $\xi$
in the Lie algebra of $H'_0$ such that $x_1=e^\xi$. Let $\nu_t = e^{-t\xi}\nu_t' = e^{-t\xi}x_t\gamma_t$.
Then $\nu_t*\nu_s=(e^{-t\xi}x_t\gamma_t)*(e^{-s\xi}\nu_s')=(e^{-t\xi}x_te^{-s\xi}\gamma_t)*\nu_s'=(e^{-t\xi}e^{-s\xi}
\nu_t')*\nu_s'=e^{-(t+s)\xi}\nu_{t+s}'=\nu_{t+s}$. This shows that $\nu_t$ is a convolution semigroup parameterized
by $\bQ_+$, which is obviously continuous. Moreover, $\nu_1=e^{-\xi}x_1\gamma_1=\gamma_1=\nu$.
By \cite[Theorem~1.5.9]{heyer}, $\nu_t$ may be extended to become a continuous convolution semigroup defined for
all $t\in\bR_+$. Because for $t\in\bQ_+^*$, $\gamma_t\in\Gamma$ and $\nu_t=e^{-t\xi}x_t\gamma_t\in\tilde{\Gamma}$, it follows
from the continuity that $\nu_t\in\tilde{\Gamma}$ for any $t\in\bR_+$. \ $\Box$
\vspace{2ex}

  The following corollary is a direct consequence of Theorem~\ref{thembeddingX} and Proposition~\ref{prmutnut}.

\begin{coro} \label{coinv}
Let $G$ be a Lie group of class $\myC$ and let $K$ be a compact subgroup. Then any $K$-right invariant
infinitely divisible distribution on $G$ is $K$-bi-invariant, and any infinitely divisible distribution on $X=G/K$
is $K$-invariant.
\end{coro}

\noindent {\bf Acknowledgement:} \ The author wishes to thank David Applebaum for some very helpful comments,
and the anonymous reviewer for the careful explanation of several points in \cite{heyer} which has led to some
important changes to the discussion on the embedding problem.

\end{document}